\documentclass{qjmam}
\startpage{1}
\yr{2018}
\vol{0}
\issue{0}

\usepackage{appendix}
\usepackage{amsmath,amssymb}
\usepackage{amsthm} 
\usepackage{chngcntr} 
\usepackage{comment}
\usepackage{float}
\usepackage[a4paper]{geometry}
\usepackage{graphicx}
\usepackage{microtype}

\usepackage[capitalize]{cleveref} 

\crefname{subsection}{subsection}{subsections} 
\crefname{app}{Appendix}{Appendices}

\DeclareMathAlphabet{\mathpzc}{OT1}{pzc}{m}{it}

\newcommand{\pd}[2]{\frac{\partial #1}{\partial #2}}

\newtheorem{proposition}{Proposition}
\newtheorem{remark}{Remark}

\everymath{\displaystyle}

\DeclareMathAlphabet\mathbit
\encodingdefault\rmdefault\bfdefault\itdefault
\DeclareOldFontCommand{\bi}{\normalfont\bfseries\itshape}{\mathbit}

\def\fakebold#1{\relax\ifvmode\leavevmode\fi%
	\ifmmode%
	\setbox0=\hbox{$#1$}%
	\else%
	\setbox0=\hbox{#1}%
	\fi%
	\kern-.02em\copy0 \kern-\wd0%
	\kern .04em\copy0 \kern-\wd0%
	\kern-.0125em\raise.02em\box0%
}%

\title[non-ideal~collapsing~cavities~and~converging~shocks]{COLLAPSING CAVITIES AND CONVERGING SHOCKS IN NON-IDEAL MATERIALS}
\author[boyd,~schmidt,~ramsey,~and~baty]{ZACHARY M. BOYD,\footnote{zachboyd@email.unc.edu}}
\def\lanl{Los Alamos National Laboratory,
New Mexico, United States}
\def\ucla{Mathematics Department, UCLA, Los Angeles, California, United States}
\def\unc{Mathematics Department, University of North Carolina, Chapel Hill, North Carolina, United States}
\def\unm{Department of Physics, New Mexico Institute of Science and Technology, Socorro, New Mexico, United States}
\address{\lanl \\
\and \unc \\
\and \ucla}
\extraauthor{EMMA M. SCHMIDT,\footnote{emschmidt@lanl.gov}}
\extraaddress{\lanl\ \and\ \unm}
\extraauthor{SCOTT D. RAMSEY,\footnote{ramsey@lanl.gov} \and ROY S. BATY\footnote{rbaty@lanl.gov}}
\extraaddress{\lanl}
\received{\recd 1 December 2018.}

\graphicspath{{./images/}} 

\begin{document}

\maketitle

\begin{abstract}
  As modern hydrodynamic codes increase in sophistication, the availability of realistic test problems becomes increasingly important. 
  In gas dynamics, one common unrealistic aspect of most test problems is the ideal gas assumption, which is unsuited to many real applications, especially those involving high pressure and speed metal deformation. Our work considers the collapsing cavity and converging shock test problems, showing to what extent the ideal gas assumption can be removed from their specification. It is found that while most materials simply do not admit simple (i.e.\ scaling) solutions in this context, there are infinite-dimensional families of materials which do admit such solutions. We characterize such materials, derive the appropriate ordinary differential equations, and analyze the associated nonlinear eigenvalue problem. It is shown that there is an inherent tension between boundedness of the solution, boundedness of its derivatives, and the entropy condition. The special case of a constant-speed cavity collapse is considered and found to be heuristically possible, contrary to common intuition. Finally, we give an example of a concrete non-ideal collapsing cavity scaling solution based on a recently proposed pseudo-Mie-Gruneisen equation of state.
\end{abstract}

\section{Introduction}

The problem of a cylindrical or spherical void converging within an ideal, inviscid fluid was first treated by Rayleigh in 1917~\cite{rayleigh} assuming incompressible flow (see also~\cite{birkhoff}). The collapsing cavity solution (and modifications thereof) has in the time since been subject to numerous additional investigations, including by Hunter~\cite{hunter} and Gilmore~\cite{gilmore} for a compressible fluid modeled by a Tait equation of state, Zwick and Plesset~\cite{zwick} with the effect of a small amount of vapor counterbalancing the collapse, and Lazarus~\cite{lazarus} for the fully compressible case in an ideal gas. The collapsing cavity solution is relevant in several practical contexts, including cavitation-induced damage, sonoluminescence~\cite{jarman1960sonoluminescence}, and as a potentially challenging test problem for the quantitative verification of inviscid compressible flow (Euler) codes.

Similarly, the problem of an infinitely strong, cylindrical or spherical shock wave converging into an ideal, compressible, inviscid fluid was first solved by Guderley~\cite{guderley} in 1942. The classical Guderley solution has subsequently been investigated by authors such as Stanyukovich~\cite{stanyukovich}, Zel'dovich and Raizer~\cite{zeldovich_and_raizer}, and Chisnell~\cite{chisnell}. Much of the attention surrounding this problem has been related to increasingly precise evaluations of its solution, its connection to potential initiating events such as curvilinear shock tubes or pistons, and applications to astrophysical processes~\cite{fink2007double}, laser fusion~\cite{atzeni2004physics,motz1979physics}, and quantitative code verification~\cite{guderley_revisited}.

\begin{figure}
	\centering
	\includegraphics[width=.7\textwidth]{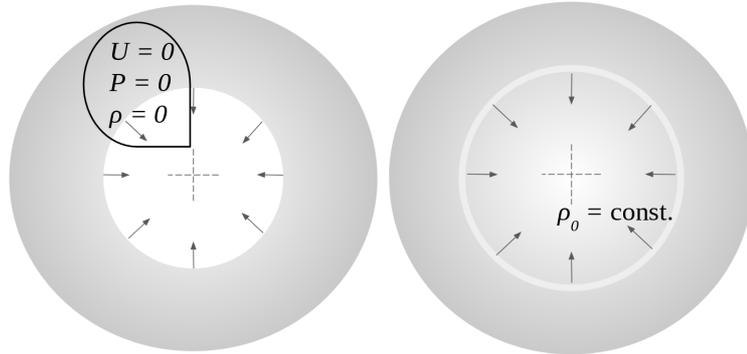}
	\caption{The Guderley (left) and collapsing cavity (right) problems. Both involve converging jumps in the velocity field. The density in the center causes the two problems to differ.}
	\label{fig}
\end{figure}

In the context of inviscid compressible flow, both of these problems feature a flow field driving a discontinuity (to avoid ambiguity, we will use the word ``jump'' to describe both processes) into a quiescent, zero-pressure material (see~\cref{fig}). 
The only physical difference between the two problems is the mass density in the undisturbed region: either zero or a non-zero constant for the collapsing cavity and converging shock problems, respectively.%
(Problems featuring power-law mass density variation in the undisturbed region have also been investigated. See, for example, Lazarus~\cite{lazarus}.)
%
With this slight difference in mind, Lazarus unified the problems into one mathematical framework though the use of a logical variable, conducted an extensive phase plane analysis of the underlying differential equations, and constructed a variety of solutions.

In his approach, Lazarus~\cite{lazarus} also showed that both problems arise from scale-invariant, self-similar transformations of the underlying governing equations and initial conditions. The converging shock and collapsing cavity problems thus represent two closely-related examples of self-similar solutions of the inviscid Euler equations. The utility of self-similar solutions has been widely disseminated, and in particular by Barenblatt:~\cite{barenblatt,barenblatt2003scaling} in addition to their well-known physical relevance (e.g., for assistance in the identification of potential scaled experiments in the contexts discussed above), they are useful as special classes of limiting behaviors known as ``intermediate asymptotics.''

Despite these advantages, existing collapsing cavity and converging shock solutions are of limited practical utility due to simplifying assumptions imposed during the formulation of their governing equations. For example, assumption of the ideal gas law is essential to the construction of these and many other scale-invariant, self-similar flows; this close correspondence owes principally to the absence of dimensional parameters within an idealized equation of state (EOS) model. However, ideal EOS laws are in general not appropriate for the characterization of important material deformation processes including tension, compression limits, non-zero low temperature sound speeds, and interatomic or intermolecular effects. These phenomena are common to real materials such as metals, granular materials, alloys, or plasmas.

The process of rigorously generalizing the EOS closure laws that allow for the existence of scale-invariant, self-similar converging shock or collapsing cavity solutions is rooted in Ovsiannikov's seminal work using symmetry analysis theory~\cite{ovsiannikov_book}; identical results and applications with enhanced discussion are provided by Holm~\cite{holm1976symmetry}, Hutchens~\cite{hutchens}, Axford~\cite{axford}, and Boyd et al.~\cite{general_euler_symmetries}. Complementary work based on physical arguments is provided by Sedov~\cite{sedov_book}, Zel'dovich and Raizer~\cite{zeldovich_and_raizer}, and Rae~\cite{rae1970analytical}. The common theme of these developments is that a necessary condition for the preservation of self-similarity (represented through either symmetry properties or dimensional arguments) is that the included EOS must assume a specific functional form that somewhat generalizes the ideal gas EOS\@. This result, however, is agnostic to any initial, boundary, or ancillary conditions that may be provided for the definition of specific flows.
 
As a result, the more subtle question regarding existence of the scale-invariant, self-similar converging shock solution for generalized EOS models was most recently addressed by Boyd et al.~\cite{guderley_general_eos}. In this work, Boyd et al.\ demonstrated that the aforementioned EOS condition was necessary but insufficient for the existence of these solutions. Further conditions related to thermodynamics and solution boundedness imposed additional constraints, so that the only possible scale-invariant, self-similar converging shock solutions with these properties are in the nonlinear eigenvalue class as discussed by Zel'dovich and Raizer~\cite{zeldovich_and_raizer}, Barenblatt~\cite{barenblatt2003scaling}, or Ramsey et al.~\cite{guderley_revisited}.
 
A key notion in the work of Boyd et al.~\cite{guderley_general_eos} is the necessity of seeking conditions for the existence of bounded, stable solutions directly from the symmetry-reduced inviscid Euler equations. This phenomenon is closely tied to existing descriptions of ``second-type similarity solutions,'' which cannot be computed based solely on dimensional or symmetry reductions. Through example, Boyd et al.\ also discuss some key features of these flows (selected from a broader list of salient properties including initial conditions, overprescribed boundary conditions, conservation principles, and the location of various characteristic curves within solution fields). While much of this phenomenology has long been understood and repeatedly applied for the case of the ideal gas, Boyd et al.\ demonstrated that a similar procedure must be followed for any generalized EOS with the correct dimensional structure arising from symmetry or other considerations.
 
Our work is thus motivated by the previous efforts concerning the existence and calculation of scale-invariant, self-similar converging shock waves. Despite the potential for a unified treatment owing to Lazarus~\cite{lazarus}, no significant corresponding attention appears to have been paid to the collapsing cavity problem in a non-ideal medium. It is thus the purpose of this work to generalize the existing results of Boyd et al.~\cite{guderley_general_eos} to a scenario that encapsulates both problems, in the style of Lazarus. To this end, we seek conditions for the existence of scale-invariant, self-similar collapsing cavity and converging shock solutions under a generalized EOS structure, with ancillary conditions including both thermodynamic soundness and boundedness everywhere within the flow field. As was the case with the converging shock solution considered in isolation, symmetry or dimensional considerations must be supplemented by a phase space singular analysis in order to arrive at a sufficiency of necessary conditions for existence of these solutions.
 
In support of our goals,~\cref{sec:governing} provides a brief review of the governing equations and conditions necessary to define the collapsing cavity and converging shock problems. A reduction of these equations to ordinary differential equations under scale-invariant, self-similar transformations is provided in~\cref{sec:sim-analysis}. \Cref{sec:odes} includes the analysis of these equations, providing the rigorous details necessary to categorize possible solutions. \Cref{sec:example} provides an example of a non-ideal collapsing cavity solution using a recently proposed pseudo-Mie-Gruneisen EOS\@. We conclude in~\cref{sec:conclusion}, with additional calculational details of potential future interest provided in two Appendices.

\section{Governing Equations and Assumptions}
\label{sec:governing}

The setting for the collapsing cavity and Guderley problems is a compressible, inviscid, adiabatic fluid. We also assume perfect cylindrical or spherical symmetry and that the fluid is at rest, i.e.\ no external forces act on it. The governing partial differential equations (PDEs) for such a flow are
\begin{align}
	\frac{\partial \rho}{\partial t} + u\frac{\partial \rho}{\partial r} + \rho\left(\frac{\partial u}{\partial r} + \frac{ku}{r}\right) = 0 \label{eqn:mass} \\
	\frac{\partial u}{\partial t} + u\frac{\partial u}{\partial r}  + \frac{1}{\rho}\frac{\partial p}{\partial r} = 0 \label{eqn:momentum}\\
	\frac{\partial p}{\partial t} + u \frac{\partial p}{\partial r} + K_S\left(\frac{\partial u}{\partial r} + \frac{ku}{r}\right) = 0 \label{eqn:energy}
\end{align}
where $\rho$, $u$, and $p$ are the density, velocity, and pressure, respectively, and $k = 1,2$ is a parameter denoting cylindrical or spherical symmetry, respectively. (See chapter one in any of~\cite{landau_and_lifschitz,courant_and_friedrichs,zeldovich_and_raizer} for derivations of~\crefrange{eqn:mass}{eqn:momentum}. Axford~\cite{axford} derives the exact form of~\cref{eqn:energy} used in this paper.) The first two equations arise from conservation of mass and momentum, and the third reflects an assumption that the material derivative of the entropy $S$ is zero. The latter assumption is not true at discontinuities; these will be handled using more general jump conditions. The material-specific properties of the flow are encoded by $K_S = K_S(p,\rho)$, which is the \emph{adiabatic bulk modulus}, defined as 
\[
	K_S = \rho \left.\pd{p}{\rho}\right|_S = \rho c^2,
\]
where $c$ is the local speed of sound. 
Using this relation, the adiabatic bulk modulus corresponding to a material can be computed from an equation of state of the form $p = p\left(\rho,S\right)$.



\Crefrange{eqn:mass}{eqn:energy} govern only smooth flows. Since we are considering shock or other discontinuous solutions, it suffices to make explicit exception for isolated jumps and require smoothness elsewhere. The corresponding jump conditions (in a zero pre-jump velocity frame) are (see, for example, Zel'dovich and Raizer~\cite{zeldovich_and_raizer})
\begin{align}
	\rho_1(u_s-u_1) = \rho_0u_s \label{eqn:jump1}\\
	\rho_0u_su_1 = p_1-p_0\\
	\rho_0u_s\left(e_1-e_0 + \frac{u_{1}^2}{2}\right) = p_1u_1\label{eqn:jump3},
\end{align}
where the subscripts $0$ and $1$ denote pre- and post-jump conditions, respectively, $u_s$ is the jump velocity, and $e$ is the specific internal energy, which can be obtained from the adiabatic bulk modulus via
\begin{equation}
	K_S \left.\pd{e}{p}\right|_{\rho} + \left.\rho\pd{e}{\rho}\right|_p = \frac{p}{\rho}
	\label{eqn:ks2e}
\end{equation}
which is derived in e.g.\ Axford~\cite{axford}. A constant of integration arising from this construction depends only on the isentrope, and is set to zero. For example, in the ideal gas law $(\gamma-1)\rho e = p$, one has $K_S = \gamma p$.

In addition to~\crefrange{eqn:jump1}{eqn:jump3}, we require an entropy condition: the jump approaches the pre-jump region supersonically but is subsonic from the perspective of the post-jump particles.
In the converging shock problem, this requirement prevents violations of the second law of thermodynamics and can be derived from the assumptions $\left.\pd{p}{\rho}\right|_S >0$, $\left.\frac{\partial^2 p}{\partial \rho^2}\right|_S>0$, and $\left.\pd{p}{S}\right|_{\rho} > 0$. (See~\cite{courant_and_friedrichs} section 65.) In the collapsing cavity problem, it is not clear that all of these assumptions hold. Our condition is, nonetheless, standard in treatments of collapsing cavity problems, e.g.\ \cite{courant_and_friedrichs} section 158 and \cite{lazarus,hunter,richtmyer-lazarus}, so we will follow this convention. We also note that other authors view this requirement as imposing a one-dimensional stability on the flow~\cite{jeffrey,burgess}.%
%

Finally, during the course of this work we will be concerned only with solutions that contain one such jump. Otherwise, we will require both smoothness and boundedness in the remainder of the solution field. In particular, $\rho$, $u$, $p$, and any derived quantities that depend on these variables (e.g., $K_S$ or the local sound speed $c$) must remain bounded as $r \to \infty$ for all $t$. 

\section{Similarity Analysis}
\label{sec:sim-analysis}

As stated in the introduction, we will be focusing on self-similar scaling solutions to the governing equations. In general, these are solutions of the form
\begin{align}
	r &= |t|^{c_0}\xi \label{eqn:scal-beg}\\
	u&= |t|^{c_1}V(\xi)\\
	\rho &=  |t|^{c_2}R(\xi)\\
	p-p_0 &= |t|^{c_3} \Pi(\xi).\label{eqn:scal-end}
\end{align}
Such solutions are called self-similar because, given the solution at some time $t = t_1$, the solution at any other time $t = t_2$ can be obtained by simply scaling the variables according to~\crefrange{eqn:scal-beg}{eqn:scal-end}. 
The essential utility of this phenomenon is reduced dimensionality in the independent variables, since the value of the solution at all points in space is needed at only a single time in order to understand the whole evolution. As we will see, self-similar processes also yield a concrete reduction to ordinary differential equations (ODEs) from the governing PDEs.

Other types of similarity solutions are defined analogously except that the solutions at distinct times are related to each other by transformations other than the self-similar scaling laws indicated by~\crefrange{eqn:scal-beg}{eqn:scal-end}.

\subsection{Self-Similar Scaling and the Bulk Modulus}

The conditions under which~\crefrange{eqn:mass}{eqn:energy} admit self-similar scaling solutions of the form indicated by~\crefrange{eqn:scal-beg}{eqn:scal-end} may be determined through various means. Of these, the simplest procedure involves direct substitution of~\crefrange{eqn:scal-beg}{eqn:scal-end} into~\crefrange{eqn:mass}{eqn:energy}, followed by selection of the constants $c_0$, $c_1$, $c_2$, and $c_3$ so that the resulting expressions collapse to ODEs in the variables $V\left(\xi\right)$, $R\left(\xi\right)$, and $\Pi\left(\xi\right)$. This is the procedure employed in the context of an ideal gas EOS, including in works by Guderley~\cite{guderley}, Stanyukovich~\cite{stanyukovich}, Zel'dovich and Raizer~\cite{zeldovich_and_raizer}, and Lazarus~\cite{lazarus}.

A complication with extending such an analysis to non-ideal materials arises in the rigorous determination of the forms of $K_S$ such that~\crefrange{eqn:mass}{eqn:energy} collapse to ODEs. This procedure is straightforward when $K_S$ assumes the ideal gas or other simple forms, but is less so when it is regarded as an arbitrary function of $\rho$ and $p$. Symmetry analysis techniques are useful in circumventing this complication.

In general, the goals of symmetry analysis of differential equations are twofold: 1) the determination of transformations under which a set of differential equations is invariant, through a local vector field or ``group generator'' representation, and 2) the use of the admissible group generator to construct a change of variables as exemplified by~\crefrange{eqn:scal-beg}{eqn:scal-end}, which are otherwise typically introduced via ansatz or other largely ad hoc considerations. For cases where the governing equations contain arbitrary functions [e.g., $K_S$ as appearing in~\cref{eqn:energy}], this framework also results in ancillary conditions these functions must satisfy to ensure the presence of various symmetries.

The symmetry analysis of~\crefrange{eqn:mass}{eqn:energy} coupled to the possible choices of $K_S\left(\rho,p\right)$ has been previously carried out by Ovsiannikov~\cite{ovsiannikov_book}, Holm~\cite{holm1976symmetry}, Hutchens~\cite{hutchens}, Axford~\cite{axford}, and Boyd et al~\cite{general_euler_symmetries}. These authors have found that~\crefrange{eqn:mass}{eqn:energy} can possess up to three independent scaling symmetries, depending on the choice of $K_S$. 
They are   
\begin{align}
  T_{1,\epsilon} &: (r,t) \rightarrow \left(e^{\epsilon}r,e^{\epsilon}t\right) \\
  T_{2,\epsilon} &: (\rho,p-p_0) \rightarrow \left(e^{\epsilon}\rho,e^{\epsilon}(p-p_0)\right) \\
  T_{3,\epsilon} &: \left(r,u,(p-p_0)\right) \rightarrow (e^{\epsilon} r, e^{\epsilon} u, e^{2\epsilon} (p-p_0)),
\end{align}
where $p_0$ is a reference pressure, $e\approx2.71$ is Euler's number and $\epsilon$ is a parameter.
Composition of these scaling operators with different values of $\epsilon$ gives new scaling operators, so that there is a three-parameter family of scaling symmetries available. It is convenient to work with local, infinitesimal transformations (called group generators) rather than global ones, which we denote by
\begin{align}
	\chi_1 &= r\pd{}{r} + t\pd{}{t} \label{eqn:X-1}\\
	\chi_2 &= \rho\pd{}{\rho} + (p - p_0)\pd{}{p} \\
	\chi_3 &= r\pd{}{r} + u\pd{}{u} + 2 (p-p_0) \pd{}{p}\label{eqn:X-3}
\end{align}
respectively. See, for example, Cantwell~\cite{cantwell} for numerous examples on how to transform between group generators and global transformations. As a brief example, one can solve the differential equation $\frac{dx}{d\tau} = \chi_1 x$ for $x=(r,t,u,\rho,p-p_0)$ and $\tau=0$ to $\epsilon$. At the end of this evolution, the $r$ and $t$ components of $x$ will have been scaled by a factor of $e^{\epsilon}$. This shows how one can link a differential operator such as $\chi_1$ to a global scaling symmetry such as $T_1$.

Linear combination of the group generators corresponds to composition of the global transformation laws. Thus, the most general symmetry we will consider is
%
\begin{equation}
	\chi = t\pd{}{t} + (1+\alpha)r\pd{}{r} + \alpha u \pd{}{u} + \beta \rho\pd{}{\rho} + (\beta + 2\alpha)(p-p_0)\pd{}{p},
	\label{eqn:generator}
\end{equation}
where we have normalized the coefficient of $t\pd{}{t}$ to 1. (See~\cref{app:A} for the case where this coefficient is $0$.)
%
\Cref{eqn:generator} represents the maximal scaling group generator admitted by~\crefrange{eqn:mass}{eqn:energy}.
Since the energy conservation relation given by~\cref{eqn:energy} contains the arbitrary function $K_S\left(\rho,p\right)$, the invariance of~\cref{eqn:energy} under the group generated by~\cref{eqn:generator} yields a conditional symmetry. In particular, the adiabatic bulk modulus must have a form compatible with the dimensional structure arising from the scaling transformations generated by~\cref{eqn:generator}. 
Invariance of~\cref{eqn:energy} under~\cref{eqn:generator} indicates $K_S$ must satisfy the PDE
\begin{equation}
	\beta\rho\pd{K_S}{\rho} + \left(2\alpha+\beta\right)\left(p-p_0\right)\pd{K_S}{p} = \left(2\alpha+\beta\right)K_S.
	\label{eqn:ks-equation}
\end{equation}
as shown in detail by Boyd et al~\cite{guderley_general_eos}. \Cref{eqn:ks-equation} indicates that for each choice of $\alpha$ and $\beta$, there corresponds some $K_S$ for which~\cref{eqn:energy} is invariant. 

Moreover, any initial or boundary conditions included in a problem formulation must also be invariant under~\cref{eqn:generator}. 
The Guderley features $\rho = \rho_0$, $u = 0$, and $p = 0$ in the undisturbed region, so that for invariance
\begin{align}
  \chi\left(\rho - \rho_0\right) = 0 &\text{ when } \rho = \rho_0 \label{eqn:iic-1}\\
  \chi\left(u - 0\right) = 0 &\text{ when } u = 0 \label{eqn:iic-2}\\
  \chi\left(p - 0\right) = 0 &\text{ when } p = 0 \label{eqn:iic-3}
\end{align}
With~\cref{eqn:generator}, these relations may be rewritten as, respectively,
\begin{align}
  \beta \rho = 0 &\text{ when } \rho = \rho_0 \label{eqn:iic-4}\\
  \alpha u = 0 &\text{ when } u = 0 \label{eqn:iic-5}\\
  (\beta + 2\alpha)(p-p_0) = 0 &\text{ when } p-p_0 = 0 \label{eqn:iic-6}
\end{align}
\Cref{eqn:iic-5} and~\cref{eqn:iic-6} are identically satisfied since $u = 0$ and $p-p_0 = 0$ in the unperturbed region, as indicated. These conditions thus provide no additional constraints on $\alpha$ and $\beta$. Conversely,~\cref{eqn:iic-4} is only satisfied for $\beta = 0$, since by definition $\rho_0 > 0$ in the unperturbed region. As a result, classical self-similar converging shock solutions must feature $\beta = 0$, and the density scaling term is eliminated from~\cref{eqn:generator}.

The collapsing cavity problem features $\rho = 0$, $u = 0$, and $p = 0$ in the undisturbed region, so that following a similar analysis as presented above, both $\alpha$ and $\beta$ are not constrained beyond any material-driven constraints imposed via~\cref{eqn:ks-equation}.
\begin{remark}
One might reasonably wonder why it is not possible to also offset $\rho$ by replacing it with $\rho-\rho_0$ in~\crefrange{eqn:X-1}{eqn:X-3}, thus relieving the restriction $\beta=0$ for the converging shock problem. While $p$ and $\rho$ do have some similarities with respect to the boundary conditions, they play fundamentally different roles in~\crefrange{eqn:mass}{eqn:energy}. In particular, whereas $p$ only appears within derivatives and in $K_S$, $\rho$ appears as a multiplicative factor, which restricts its ability to be offset additively without altering the PDEs themselves. This is confirmed rigorously in, for example,~\cite{guderley_general_eos}.
\end{remark}
The possible values of $K_S$ for both problems are summarized in~\cref{tab:eos-restrictions}.
\begin{table}[H]
	\centering
	\begin{tabular}{l|l|l}
		$K_S$							&	Collapsing cavity				&	Guderley \\
		\hline
		$(p-p_0)f(p,\rho)$					&	$\alpha=\beta=0$		&	$\alpha=\beta=0$	\\
		$(p-p_0)f\left( (p-p_0)\rho^{-\lambda} \right)$	&	$\beta+2\alpha = \lambda\beta$ 	&	$\alpha=\beta=0$	\\
		$(p-p_0)f(\rho)$					&	$\beta=0$			&	$\beta=0$		\\
		$(p-p_0)\gamma$					&	No constraints on $\alpha$ or $\beta$				&	$\beta=0$
	\end{tabular}
	\caption{The correspondence between values of $K_S$ and choices of $\alpha$ and $\beta$ in~\cref{eqn:generator} for the two problems under consideration. Here, $f$ is an arbitrary smooth function, $\lambda \in \mathbb{R}$, and $\gamma$ is the adiabatic index from the ideal gas law. The difference between the two columns arises because~\crefrange{eqn:iic-3}{eqn:iic-6} restrict the Guderley problem but not the collapsing cavity problem.}
	\label{tab:eos-restrictions}
\end{table}

We summarize these facts in the following proposition.
\begin{proposition}
	The scaling solutions of the collapsing cavity and converging shock problems are invariant under the operator~\cref{eqn:generator} and satisfy the restrictions given in~\cref{tab:eos-restrictions}. As will be seen later, those solutions which do exist satisfy a nonlinear eigenvalue problem.
\end{proposition}

\Cref{tab:eos-restrictions} may be physically interpreted as follows: the three scaling groups in~\crefrange{eqn:X-1}{eqn:X-3} represent three ways to reduce the number of variables in the problem. Normalizing to the lead coefficient to $1$ in~\cref{eqn:generator} used up one (redundant) degree of freedom. It will be seen later that one of the two remaining degrees of freedom is needed for the solution of a nonlinear eigenvalue problem, so only one degree of symmetry remains to be spent.

One place where the symmetry can be broken is in the form of $K_S$, since heuristically, we lose one degree of scaling for each dimensional parameter appearing in it (since dimensional parameters induce characteristic scales, thus removing scale freedom). The units of $K_S$ are the same as those of pressure, so that the ideal gas law, for instance, does not require the introduction of any dimensional parameters. In contrast, the second and third rows of~\cref{tab:eos-restrictions} need a dimensional parameter so that the output of $f$ can be rendered dimensionless. The first row in general needs two dimensional parameters to cancel the dimensions of its arguments, which is why it has comparatively little scaling symmetry.

The other way that scale invariance is lost is through initial or boundary conditions, which can also contain dimensional parameters. This is not relevant to the cavity problem, but the constant undisturbed density condition in the converging shock problem does restrict us to symmetries in which there is no density scaling (i.e., $\beta = 0$). One would think that this would make the converging shock problem impossible to solve in the case where $K_S$ also includes a dimensional parameter, but there is one important exception. When the dimensional constants arising from the initial or boundary conditions and bulk modulus constraint are of the same kind, the case considered in the third row of~\cref{tab:eos-restrictions} arises. Other than that (and the fourth row, which is a special case of the third), there are indeed too many degrees of freedom used up, and the converging shock problem cannot be solved by similarity methods in these cases. While $\alpha=\beta=0$ does give a reduction to similarity variables, this corresponds to a constant jump speed, and it was shown in~\cite{guderley_general_eos} that there is no constant shock speed converging shock solution.

Of course, while~\cref{tab:eos-restrictions} describes the classes of bulk moduli that allow scaling solutions respecting the governing PDEs, jump conditions, and boundary conditions of the collapsing cavity and converging shock problems, further discrimination is necessary to identify specific forms of $K_S$ satisfying fundamental thermodynamic criteria (such as having positive sound speed) and fitting empirical data. A full exploration of these issues is outside the scope of the present study, although we note the (non-obvious) fact that there do exist forms of $K_S$ which match the forms in~\cref{tab:eos-restrictions}, are thermodynamically reasonable in the required ranges, are physically relevant, and differ substantially from the ideal gas law (see~\cite{frank}). A useful guide to some of the issues to check for and how to correct them is~\cite{segletes}.

In any event, the similarity variables associated with~\cref{eqn:generator} 
may be determined by 
finding the invariant functions of the group generator. In particular, for an arbitrary function $F\left(r,t,u,\rho,p\right)$, the invariant function condition given by
\begin{equation}
  \chi F\left(r,t,u,\rho,p\right) = 0, \label{eqn:invfunc}
\end{equation}
gives rise to the characteristic system
\begin{equation}
  \frac{dt}{t} 
  = \frac{dr}{\left(1+\alpha\right)r} 
  = \frac{du}{\alpha u} 
  = \frac{d\rho}{\beta \rho} 
  = \frac{dp}{\left(\beta + 2\alpha\right)\left(p-p_0\right)} \label{eqn:chareqs}
\end{equation}
The constants of integration arising from the solution of~\cref{eqn:chareqs}
are invariant functions of the group generated by~\cref{eqn:generator}, and are interpreted as similarity variables through which to transform~\crefrange{eqn:mass}{eqn:energy} to ODEs. Solving~\cref{eqn:chareqs} thus yields 
\begin{align}
	r &= t^{\alpha+1}	\xi\label{eqn:xi}\\
	\rho &=  t^{\beta}	R(\xi)\label{eqn:R}\\
	u&= t^{\alpha}	V(\xi) \label{eqn:V} \\
	p-p_0 &= t^{2\alpha+\beta}\Pi(\xi) \label{eqn:Pi}
\end{align}
where the invariant functions $\xi$, $V$, $R$, and $\Pi$ are the constants of integration. 
In assuming these similarity forms,~\crefrange{eqn:mass}{eqn:energy} reduce to
\begin{align}
	0 &= -\beta R + (\alpha + 1)R'\xi + VR' + R\left(V' + k\frac{V}{\xi}\right) \label{eqn:implicit-odes-beg} \\ 
	0 &= -\alpha V + (\alpha + 1)V'\xi + VV' +\frac{\Pi'}{R} \\
	0 &= -(\beta + 2\alpha)\Pi + (\alpha + 1)\Pi'\xi + V\Pi' + f \Pi\left(V' + k\frac{V}{\xi}\right) \label{eqn:implicit-odes-end}
\end{align}
where $f = K_S/(p-p_0)$.

\subsection{Jump Conditions}
\label{sec:jump}

While Ovsiannikov~\cite{ovsiannikov_book} and Boyd et al.~\cite{general_euler_symmetries} conduct an analysis that is restricted to the case of smooth solutions of~\crefrange{eqn:mass}{eqn:energy}, it turns out that in practice, the jump conditions usually have the same similarity behavior as the associated PDEs.
This is also true in our particular case.

In the case of a collapsing cavity,~\crefrange{eqn:jump1}{eqn:jump3} reduce to

\begin{align}
	\rho_1(u_s-u_1) = 0 \label{eqn:cavjump1} \\
	0 = p_1 \label{eqn:cavjump2} \\
	0 = p_1u_1 \label{eqn:cavjump3} ,
\end{align}
which yields $u_1=u_s$ and $p_1=0$ as the post-jump conditions. There are no restrictions on the post-jump density, except those imposed by the scaling assumption. 
The scaling behavior of these jump conditions is determined by substituting~\crefrange{eqn:R}{eqn:Pi} into~\crefrange{eqn:cavjump1}{eqn:cavjump3} to yield
\begin{align}
	|t|^{\beta} R_1(u_s - |t|^{\alpha} V_1) = 0 \label{eqn:scavjump1} \\
	0 = |t|^{2\alpha+\beta}\Pi_1 \label{eqn:scavjump2} \\
	0 = |t|^{3\alpha+\beta}\Pi_1 V_1 \label{eqn:scavjump3} ,
\end{align}
For~\crefrange{eqn:scavjump1}{eqn:scavjump3} to reduce to a statement in terms of the similarity variables, the jump speed $u_s$ must obey
\begin{equation}
  u_s \propto |t|^{\alpha} \label{eqn:jumpvel}
\end{equation}
such that the jump trajectory is given by
\begin{equation}
  r_s = \xi_s |t|^{\alpha+1} \label{eqn:traj}
\end{equation}
where $\xi_s$ has been judiciously selected so as to absorb both proportionality and integration constants. \cref{eqn:traj} also indicates that the jump stays at the fixed coordinate $\xi_s$ in similarity space.

In the case of the converging shock problem, the presence of the energy terms complicates matters somewhat, and we cannot fully compute the post-jump conditions without knowing more about the explicit form of $K_S.$ Nevertheless, we can confirm that~\crefrange{eqn:jump1}{eqn:jump3} are invariant under the scaling induced by $\chi$. To see how $e$ scales, we consider~\cref{eqn:ks2e} and apply~\crefrange{eqn:xi}{eqn:Pi}, yielding
\begin{equation}
	\Pi f \pd{e}{\Pi} + R\pd{e}{R} = |t|^{2\alpha} \frac{\Pi}{R}.
\end{equation}
Therefore, we find that $E = |t|^{-2\alpha}e$ is a function $\xi$ alone. With this fact in hand, one easily substitutes~\crefrange{eqn:R}{eqn:Pi} into~\crefrange{eqn:jump1}{eqn:jump3}, replacing $e$ with $t^{2\alpha}E$, and finds 
\begin{align}
	|t|^{\beta} R_1 \left(u_s - |t|^{\alpha} V_1\right) &= |t|^{\beta} R_0 u_s \label{eqn:simjump1}\\
	|t|^{\alpha+\beta} R_0 u_s V_1 &= |t|^{2\alpha+\beta}\Pi_1 -|t|^{2\alpha+\beta} \Pi_0\\
	|t|^{\beta} R_0 u_s 
	\left(|t|^{2\alpha}\left(E_1-E_0\right) + |t|^{2\alpha} \frac{V_1^2}{2}\right) 
	&= \Pi_1 V_1 |t|^{3\alpha+\beta} \label{eqn:simjump3},
\end{align}
or
\begin{align}
R_1(V_s-V_1) &= R_0 V_s\\
R_0 V_s V_1 &= \Pi_1-\Pi_0\\
R_0 V_s\left(E_1-E_0 + \frac{V_1^{2}}{2}\right) &= \Pi_1 V_1,
\end{align}
where the scaled shock speed $V_s \propto u_s |t|^{-\alpha}$ from~\cref{eqn:jumpvel}, and~\cref{eqn:traj} holds so that the jump conditions reduce to a statement in terms of similarity variables. 

While we cannot solve the jump conditions for completely unspecified $K_S$, they are at least well-defined in the similarity space. The main use of the jump conditions in our analysis is simply to provide an initial condition for the integration of~\crefrange{eqn:implicit-odes-beg}{eqn:implicit-odes-end}, so well-definedness is sufficient for our purposes.

\section{Analysis of ODEs}
\label{sec:odes}

Starting from~\crefrange{eqn:implicit-odes-beg}{eqn:implicit-odes-end}, we can isolate the derivatives. To do this, we observe that~\crefrange{eqn:implicit-odes-beg}{eqn:implicit-odes-end} are linear in the derivatives, so we can write the equivalent matrix equation
\begin{equation}
	\left( 
	\begin{array}{ccc}
		X & R		& 0		\\
		0 & X		& R^{-1}	\\
		0 & \Pi f	& X	
	\end{array}
	\right)
	\left(
	\begin{array}{c}
		R'\\
		V'\\
		\Pi'
	\end{array}
	\right)
	=
	\left(
	\begin{array}{c}
		\beta R - k \frac{RV}{\xi}\\
		\alpha V\\
		\left( \beta+2\alpha \right) \Pi - k \frac{V}{\xi}\Pi f
	\end{array}
	\right)
\end{equation}
where $X = (\alpha + 1) \xi + V$ is the group velocity in the reference frame of the particle at $\xi$. Inverting the matrix then yields
\begin{equation}
	\frac{d}{d\xi}
	\left(
	\begin{array}{c}
		R \\
		V \\
		\Pi
	\end{array}
	\right)
	=
	\frac{\left( \beta+2\alpha \right)\Pi - RV\left(\alpha X + k\frac{C^2}{\xi}\right)}{RX\left(X^2 - C^2\right)}
	\left(
	\begin{array}{ccc}
		R\\
		-X\\
		RX^2
	\end{array}
	\right)
	+
	\left(
	\begin{array}{c}
		\left(\beta - k\frac{V}{\xi}\right)\frac{R}{X} \\
		0 \\
		\alpha RV
	\end{array}
	\right)
	\label{eqn:ODEs}
\end{equation}
where $R C^2 = \Pi f$ gives the scaled sound speed. The form of the ODEs given in~\cref{eqn:ODEs} is particularly conducive to the treatment of global existence, but there are other forms that yield different insights (see e.g.~\cref{app:conservation-form}).

Consistent with the definitions of the collapsing cavity and converging shock problems disseminated in the existing literature (and summarized in~\cref{sec:governing}), we will seek bounded, global, smooth solutions of~\cref{eqn:ODEs}.
\begin{remark}
Most authors consider only the smooth solutions. Lazarus~\cite{lazarus} provides an in-depth investigation into the non-smooth case, and~\cite{lazarus-stability} found that only the smooth solutions are stable.
\end{remark}
However, we will show that in general there exist two points within any such solution: at one of which the denominator of~\cref{eqn:ODEs} acquires a positive value and at the other of which it takes a negative value. The intermediate value theorem then implies that the denominator vanishes somewhere between those two points, potentially ruling out solutions of the form we seek.

The first point is one very far from the origin. As $r\to\infty$, by their definitions $\xi\to\infty$ and thus $X\to\infty$ (since $V$ is bounded as $u$ is assumed bounded). By the same assumption, and by its definition, $C$ also stays bounded as $r\to\infty$. So, far from the origin, the denominator in~\cref{eqn:ODEs} is positive. 

The second point is immediately behind the jump (i.e., at $r = r_s$ or $\xi = \xi_s$). Stability (see Sec.~\ref{sec:governing}) implies that at this point the jump will propagate with a velocity less than the (physical) sound speed $c$, in a frame of reference where the particles are stationary, i.e. 
\begin{equation}
	|u_s - u| < |c|.
	\label{eqn:stability}
\end{equation} 
From the developments of~\cref{sec:jump}, the jump stays at a fixed $\xi$ coordinate, $\xi_s$. Thus we use~\cref{eqn:traj} to compute 
\begin{equation}
	u_s = \frac{dr_s}{dt} = -(\alpha+1)\xi_s |t|^\alpha
	\label{eqn:shock-speed}
\end{equation} 
Using this relation,~\cref{eqn:V}, and the relation between the physical and scaled sound speeds $c = |t|^{\alpha} C$,~\cref{eqn:stability} becomes
\begin{equation}
	\left|-(\alpha+1)\xi_s |t|^{\alpha} - V |t|^{\alpha}\right| < C |t|^\alpha, \label{eqn:simstability}
\end{equation} 
or
\begin{equation}
	\left|X \right| < C , \label{eqn:simstability2}
\end{equation}
using the definition of $X$ given above. \Cref{eqn:simstability2} indicates that $|X| < C$ at the position of the jump, $\xi = \xi_s$. The denominator in~\cref{eqn:ODEs} is therefore negative immediately behind the jump. 
Finally, the intermediate value theorem implies that the denominator in~\cref{eqn:ODEs} must vanish somewhere for all smooth solutions. 

Our attention next turns to the question of whether it is possible for the numerator in~\cref{eqn:ODEs} to vanish simultaneously with the denominator, which is the only remaining hope for a smooth, bounded solution to~\cref{eqn:ODEs}. Specifically we want to satisfy
\begin{equation}
  (2\alpha + \beta)\Pi - RV\left( \alpha X + k \frac{C^2}{\xi} \right) = 0 \text{ whenever $X^2=C^2$.}
	\label{eqn:num-cond}
\end{equation}
\Cref{eqn:num-cond} is the nonlinear eigenvalue problem for the collapsing cavity and converging shock problems. It is a single constraint with two free parameters, $\alpha$ and $\beta$. Ostensibly, at least one of these parameters may be determined so that given a solution for $R$, $V$, and $\Pi$,~\cref{eqn:num-cond} may be satisfied simultaneously with the vanishing denominator condition $X^2 = C^2$. 
Conversely, if neither $\alpha$ nor $\beta$ is free to be chosen, then in general~\cref{eqn:num-cond} cannot be satisfied simultaneously with $X^2 = C^2$, and no smooth, bounded solution can exist. 

Depending on the included material model represented by the choice of $K_S$, constraints on $\alpha$ and $\beta$ are summarized in~\cref{tab:eos-restrictions}. By comparing variables and constraints, one expects there to be no smooth, bounded solution to the collapsing cavity and converging shock problems except for those cases where~\cref{tab:eos-restrictions} imposes at most one constraint on $\alpha$ and $\beta$. In particular, 
\begin{itemize}
\item For the collapsing cavity problem, we expect solutions to exist in the cases corresponding to the last three rows of~\cref{tab:eos-restrictions}.
\item For the converging shock problem, we expect solutions to exist in the cases corresponding to the last two rows of~\cref{tab:eos-restrictions}.
\end{itemize} 
Conversely, no solutions are expected to exist for either problem in any case where $\alpha = \beta = 0$, since there is no free parameter to tune. This result includes that previously given for the converging shock problem in~\cite{guderley_general_eos}, but has now been expanded to include collapsing cavity scenarios. The case where $\alpha = \beta = 0$ is especially notable considering that it is the only case which applies to any choice of $K_S$. This result thus definitely confines the existence of any collapsing cavity or converging shock solutions (with the desired properties) to the specialized classes of materials indicated in the second, third, and fourth rows of~\cref{tab:eos-restrictions}.

Beyond this general result,~\cref{tab:eos-restrictions} includes at least two special cases. The first of these is the collapsing cavity problem in an ideal gas (the fourth row of~\cref{tab:eos-restrictions}), where there are no constraints except for~\cref{eqn:num-cond}; i.e., both $\alpha$ and $\beta$ are unconstrained.  Perhaps the most complete investigation of this problem owes to Lazarus~\cite{lazarus}. In this work, Lazarus imposes an additional constraint on the entropy of the flow, thus obtaining a class of nonlinear eigenvalue solutions similar to the converging shock solutions. However, based on our work, it is possible that the ideal gas collapsing cavity problem has two-parameter families of solutions for each choice of the adiabatic index $\gamma$. This possibility remains a potentially interesting avenue for future research.

A second special case includes any scenario where $\alpha = 0$ but $\beta$ is nonzero. While these scenarios are definitely ruled out for converging shock problems, their existence is a possibility for collapsing cavity problems (e.g., where an inherent velocity scale is somehow included in the problem formulation). Any solutions in this class yield a constant velocity jump, as has been physically observed in~\cite{constant-velocity}. Solutions of this type have been noted by Lazarus~\cite{lazarus-stability} and Thomas et al.~\cite{constant-velocity}. However, the explicit construction and utilization of solutions of this type appear to represent yet another avenue of future research potentially inspired by our work.

\section{Example Solution}
\label{sec:example}

\begin{figure}
  \centering
  \includegraphics[width=3in]{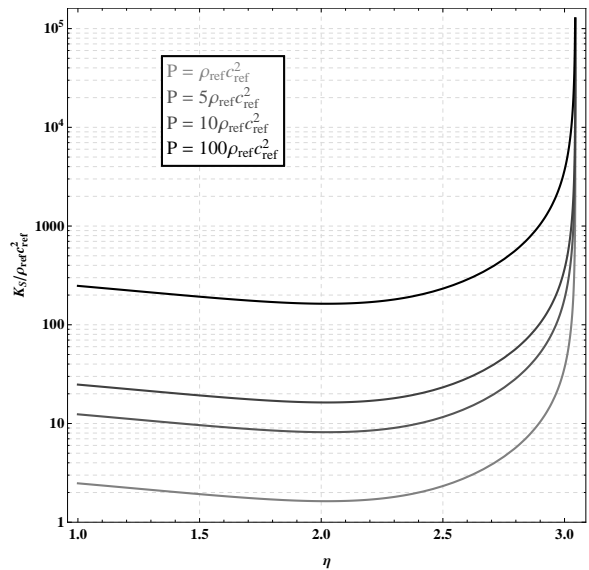}
  \caption{Normalized bulk modulus versus normalized density at different pressures $P$ for the EOS used in~\cref{sec:example}. In this formulation, $\rho_{\mathrm{ref}}$ and $c_{\mathrm{ref}}$ are the reference density and sound speed corresponding to a linear shock speed-particle speed construction of the Mie-Gruneisen EOS; these parameters may be used to select a material-specific scale for the bulk modulus, and are entirely independent of the additional EOS parameters $q$ and $s$.}
  \label{ffig3}
\end{figure}

As a concrete example of the previous developments, consider the collapsing cavity scenario in the context of the pseudo-Mie-Gruneisen EOS defined by Ramsey et al.~\cite{frank}. The bulk modulus corresponding to this EOS (see~\cref{ffig3}) is given by
\begin{equation}
  K_S = p c_1 \left( c_2 + \frac{\left(\eta - c_3\right)^2}{\eta_{\mathrm{max}} - \eta} \right)
  \label{eqn:frank-KS}
\end{equation}
where $\eta = \rho/\rho_{\mathrm{ref}}$ for some reference density $\rho_{\mathrm{ref}}$, and $c_1$, $c_2,$ and $c_3$ are constants given by
\begin{align}
  c_1 = \frac{1-4s}{4q(q-2)(s-1)} \\
  c_2 = q\left( 4(s-1) - q (2s-1) \right) \\
  c_3 = \frac{q+s-1}{s-1}
\end{align}
for parameters $s\in [0,\infty)$ and $q\in [0,1]$. Here, $s$ corresponds to the slope in a linear shock-speed, particle-speed Hugoniot relation. The parameter $q$ is a tuning constant indicative of the fractional position of a maximum in $K_S$ for the Mie-Gruneisen EOS. See~\cite{frank} for more details on both of these parameters. The value of $\eta_{\mathrm{max}}$ is then given by $\frac{s}{s-1}$. As extensively detailed by Ramsey et al.~\cite{frank}, this bulk modulus was designed to capture several key physical features of the generally applicable Mie-Gruneisen EOS, while still taking a form so as to allow for the construction of various similarity solutions (when coupled to the inviscid Euler equations). Indeed, in the context of this work, it matches the conditions in~\cref{tab:eos-restrictions} with $\beta=0$. 
We therefore consider the nonlinear eigenvalue solution of the cavity problem with $s=1.489$ and $q=0.25$ (as given by Ramsey et al.~\cite{frank} for a parameterization of copper). We also take $\rho_{\mathrm{ref}}=1$ so that all computed results for the density may be rescaled by any arbitrary $\rho_{\mathrm{ref}}$.

With~\cref{eqn:frank-KS}, \cref{eqn:ODEs} can be solved numerically given the appropriate initial conditions. As given by~\crefrange{eqn:cavjump1}{eqn:cavjump3}, for the case of a collapsing cavity the post-jump conditions are $u_1 = u_s$ (the free surface velocity is the flow velocity evaluated at the location of the free surface) and $p_1 = 0$ (the pressure at the surface of the cavity is zero); the density on the cavity surface is not constrained and will be taken as the reference value $\rho_{\mathrm{ref}}=1$. Transforming to similarity variables using~\crefrange{eqn:xi}{eqn:Pi}, the initial conditions are:
\begin{align}
  R_1 = R\left(\xi = \xi_s \right) = 1 \label{eqn:cavjumpsim1} \\
  V_1 = V\left(\xi = \xi_s \right) = -1 \label{eqn:cavjumpsim2} \\
  \Pi_1 = \Pi\left(\xi = \xi_s \right) = 0 \label{eqn:cavjumpsim3}
\end{align}

\begin{figure}
  \centering
  \includegraphics[width=4in]{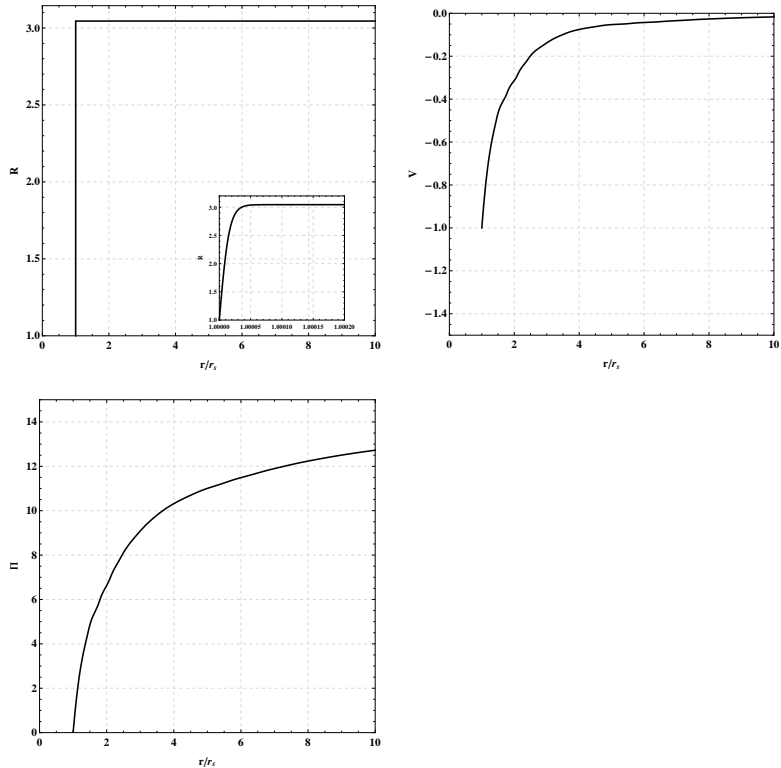}
  \caption{Scaled density ($R$), velocity ($V$), and pressure ($\Pi$) for the collapsing cavity problem for the non-ideal equation of state described in~\cref{sec:example}. The density rises very sharply but is, in fact, continuous. These curves were obtained using Wolfram Mathematica 9~\cite{Mathematica}, with $\alpha=-0.852896$.}
  \label{ex-sol}
\end{figure}

With the EOS parameterization given by~\cref{eqn:frank-KS} and the collapsing cavity initial conditions given by~\crefrange{eqn:cavjumpsim1}{eqn:cavjumpsim3}, \cref{eqn:ODEs} may be solved using any of a variety of numerical integration packages for ODEs. Using Wolfram Mathematica~9~\cite{Mathematica}, it is found via a shooting method that $\alpha = -0.852896$
satisfies the aforementioned nonlinear eigenvalue problem, and enables a smooth solution of~\cref{eqn:ODEs}. The similarity variables $R$, $V$, and $\Pi$ corresponding to this solution are depicted in~\cref{ex-sol} 

The solution depicted in~\cref{ex-sol} includes several notable features. 
Perhaps foremost is the behavior of the density field, which increases sharply (but continuously) from the interface value of $R = 1$ to the isentropic compression limit $R = \frac{s}{s-1}$ (corresponding to physical values $\rho_{\mathrm{ref}}$ at the cavity surface and the limiting value $\rho_{\mathrm{ref}}\frac{s}{s-1}$). The velocity and pressure fields vary more gradually and monotonically from their interface values. For the pressure field, this behavior is qualitatively different than that typically observed for ideal gas solutions (see, for example, Zel'dovich and Raizer~\cite{zeldovich_and_raizer}, Lazarus~\cite{lazarus}, or Ramsey et al.~\cite{guderley_revisited,frank}); however, following from Yousaf's~\cite{yousaf1986imploding} definitive resolution of the Fujimoto and Mishkin and Lazarus debate~\cite{lazarus,fujimoto1978analysis,mishkin1978analysis,lazarus1980comments,mishkin1980reply}, the existence of a universal pressure maximum is not necessarily to be expected. The flow field behaviors depicted in~\cref{ex-sol} are thus physically plausible. 

Finally, it is worth noting that the solution depicted in~\cref{ex-sol} is qualitatively different than the ideal gas collapsing cavity solutions disseminated by Lazarus~\cite{lazarus}, owing principally to our choice of $\beta = 0$. Lazarus~\cite{lazarus} instead requires 
\begin{align}
  \beta = \frac{-2 \left(\frac{1}{\alpha+1} - 1\right)}{\gamma - 1}
\end{align}
of his collapsing cavity solutions (where $\gamma$ is the ideal gas adiabatic index), which he notes corresponds to constant entropy behind the cavity surface. Our dispensing with this requirement no doubt has implications for the behavior of the flow field, including some of the notable features discussed above. Nonetheless our solution also appears to be somewhat more general, in that it features one less constraint on the flow field.

\section{Conclusion}
\label{sec:conclusion}
 
We have developed a comprehensive classification of those forms of the adiabatic bulk modulus $K_S$ which enable exact self-similar, scale-invariant solutions to the collapsing cavity and converging shock problems. These correspond to classes of generally non-ideal materials for which exact hydrodynamic self-similar scaling is expected to occur, at least in a certain physical regime. There are eight qualitatively-different cases, which are summarized in~\cref{tab:eos-restrictions}. We found that, while self-similar scaling solutions are not expected to exist for an arbitrarily-chosen adiabatic bulk modulus, there is an infinite-dimensional family of adiabatic bulk moduli for which self-similar scaling will occur, thus showing that such phenomena are not limited to the ideal gas case, as is commonly supposed. A recent application of this fact to the converging shock problem is examined in~\cite{frank}, and a complementary example in the context of the collapsing cavity problem is provided in this work.

Moreover, we also showed that some of the features of the corresponding ideal gas problems are common to all such problems, including the basic form of the resulting ODEs, the need to search for a nonlinear eigenvalue, and the tension between bounded state variables, bounded derivatives, and the entropy condition. This tension is at the root of the need for nonlinear eigenvalue analysis, and is a calling card of ``second-type similarity solutions,'' as discussed extensively by Zel'dovich and Raizer~\cite{zeldovich_and_raizer} and Barenblatt~\cite{barenblatt,barenblatt2003scaling}.

The results of our work can also be compared with those of Lazarus~\cite{lazarus,lazarus-stability} and Thomas et al.~\cite{constant-velocity}, who provide extensive numerical studies of the collapsing cavity and converging shock problems in conjunction with the ideal gas EOS model. Lazarus in particular identifies various flows that are not expected to exist based on stability arguments. Lazarus' conclusions are based primarily on properties of the linearized flow equations (written in terms of similarity variables), and are thus distinct from the thermodynamic jump condition used in our work. As such, our work is not inconsistent with these previous developments, and should be viewed as a minimum set of requirements for existence of collapsing cavity and converging shock self-similar scaling solutions. However, as demonstrated in previous work, additional criteria may exist that rule out existence of certain solutions. Examination of these criteria in the context of a general EOS closure law represents yet another potential extension of our work.

Finally, the numerical example provided in Sec.~\ref{sec:example} raises a number of interesting questions that could form the basis of future investigations. The adiabatic bulk modulus chosen in this example enables direct comparison between our results and those of Ramsey et al.~\cite{frank} for the converging shock problem. Direct comparison to the Lazarus~\cite{lazarus} results for a collapsing cavity in an ideal gas are complicated by an inconsistency in the assumed value of the density scaling constant $\beta$, as discussed in Sec.~\ref{sec:example}. Because of this inconsistency our solution does not limit to Lazarus' - and, in fact, reflects different flow physics - but additional examples could be easily devised to explore both this potential connection and the other canonical cases appearing in~\cref{tab:eos-restrictions}.

\section*{Acknowledgements}
This work was performed under the auspices of the United States Department of Energy by Los Alamos National Security, LLC, at Los Alamos National Laboratory under contract DE-AC52-06NA25396. Z.\ Boyd was additionally supported by the Department of Defense (DoD) through the National Defense Science \& Engineering Graduate Fellowship (NDSEG) Program and Army Research Office awards W911NF-18-1-0244 and W911NF-16-1-0356. Further funding for Z.\ Boyd was supplied by the James S.\ McDonnell Foundation 21st Century Science Initiative (grant \#22020315).
The authors thank E.\ J.\ Albright and J.\ Schmidt for valuable insights on these topics, as well as anonymous reviewers, who helped especially with~\cref{app:conservation-form}.

\begin{appendices}
	\crefalias{section}{app}

	\section{The case \texorpdfstring{$a_1 = 0$}{a1=0}}
	\label{app:A}
	In~\cref{eqn:generator}, we assumed $a_1 \ne 0$, so we now briefly consider the opposite case. A general linear combination of $\chi_1$, $\chi_2$, and $\chi_3$ has the form
\begin{equation}
	\chi = a_1 t\pd{}{t} + (a_1+a_3)r\pd{}{r} + a_3 \pd{}{u} + a_2 \rho\pd{}{\rho} + (a_2 + 2a_3)(p-p_0)\pd{}{p},
	\label{eqn:generator-full}
\end{equation}
	Assuming $a_1=0$ and $a_3\ne 0$, we get the similarity variables 
	\begin{align*}
		t &= t\\
		u &= rV(t)\\
		\rho &= r^{a_2/a_3} R(t)\\
		p-p_0 &= r^2 \Pi(t).
	\end{align*}
	Since such a choice of variables will not admit a solution of the kind we are seeking (for example, the velocity is linear in space and therefore cannot be zero before the jump), we can safely dismiss this case. Alternatively, when $a_1=a_3=0$, there is no scaling at all on $r$ or $t$, so that a reduction to ODEs does not occur. Therefore, we may reject this case as well.

	\section{A ``conservation'' form of the ODEs}
	\label{app:conservation-form}

	It is worth noting that~\cref{eqn:ODEs} can be reformulated as
	\begin{align}
	  (RX)' &= \left((1+\alpha)(1+k) + \beta - k \frac{X}{\xi}\right) R \label{c1}\\
	  \left( RX^2+\Pi \right)' &= \left(2 + k + (3+k)\alpha + \beta - k \frac{X}{\xi} - \alpha(1+\alpha)\frac{\xi}{X}\right) RX \label{c2}\\
	  R'C^2 - \Pi' &= \beta R \frac{C^2}{X} - (\beta+2\alpha) \frac{\Pi}{X} \label{c3}
	\end{align}
	\Cref{c1}~is a statement regarding a mass flux, and \cref{c2}~is the corresponding flux law for momentum. (To see this, integrate the equations over some range of $\xi$.) \Cref{c3} can be related to entropy by considering 
	\[
	  |t|^{2\alpha+\beta} (C^2R' - \Pi') = \left.\pd{p}{\rho}\right|_S \frac{d\rho}{d\xi} - \frac{dp}{d\xi}  = -\left.\pd{p}{S}\right|_\rho \frac{dS}{d\xi},
	\]
	where the last equality follows from the identity 
	\[\frac{dp}{d\xi}  = \left.\pd{p}{S}\right|_\rho \frac{dS}{d\xi} + \left.\pd{p}{\rho}\right|_S \frac{d\rho}{d\xi}.\]
	  A potentially fruitful avenue of future research may include connecting the structure of~\crefrange{c1}{c3} to our earlier conclusions regarding boundedness and thermodynamic soundness.

\end{appendices}

\nocite{*}
\bibliographystyle{my-qjmam}
\bibliography{cavity}

\end{document}